\documentclass[a4paper, 12pt]{article}

\usepackage{graphics}
\usepackage{graphicx}
\usepackage{amsmath}
\usepackage{amsfonts}
\usepackage{indentfirst}

\usepackage{fancyhdr} 
\fancypagestyle{plain}

\begin{document}
\date\empty
\title{Numerical solutions to boundary value problem for anomalous
       diffusion equation \\ with Riesz-Feller fractional operator}
\author{ Mariusz Ciesielski, Jacek Leszczynski}
\maketitle

\begin{center}
Institute of Mathematics and Computer Science, \\[0pt]
Czestochowa University of Technology \\[0pt]
ul. Dabrowskiego 73, 42-200 Czestochowa \\[0pt]
e-mail: mariusz@imi.pcz.pl, jale@imi.pcz.pl
\end{center}

\noindent 
In this paper, we present a numerical solution to an ordinary differential 
equation of a fractional order in one-dimensional space. The solution to 
this equation can describe a~steady state of the process of anomalous diffusion.
The process arises from interactions within complex and non-homogeneous
background. We present a numerical method which is based on the finite
differences method. We consider a~boundary value problem (Dirichlet
conditions) for an equation with the Riesz-Feller fractional derivative.
In the final part of this paper, same simulation results  are shown.
We present an example of non-linear temperature profiles in nanotubes
which can be approximated by a solution to the fractional differential equation.


\section{Introduction}

\setcounter{equation}{0} 
Anomalous diffusion is a phenomenon strongly connected with interactions 
within complex and non-homogeneous background. This phenomenon is observed 
in transport of a fluid in porous materials, in chaotic heat baths, 
amorphous semiconductors, particle dynamics inside a polymer network,
two-dimensional rotating flow and also in econophysics. The phenomenon 
of anomalous diffusion deviates from the standard diffusion behaviour. 
In opposite to the standard diffusion where a~linear form in
the mean square displacement $\left\langle \,x^{2}\left( t\right)
\right\rangle \sim k_{1}t $ of the diffusing particle over time occurs,
the anomalous diffusion is characterized by a non-linear one $\left\langle
\,x^{2}\left( t\right) \right\rangle \sim k_{\gamma }t^{\gamma} $, for $%
\gamma \in (0,2]$. In this phenomenon, there may exist a dependence $\left\langle
\,x^{2}\left( t\right) \right\rangle \rightarrow \infty $, which is
characterized by occurrence of rare but extremely large jumps of diffusing
particles -- well-known as the Levy motion or the Levy flights. An ordinary
diffusion follows Gaussian statistics and Fick's second law for finding
the running process at time $t$ whereas anomalous diffusion follows non-Gaussian
statistic or can be interpreted as the Levy stable densities.

Many authors proposed models which base on linear and non-linear forms of
differential equations. Such models can simulate anomalous diffusion but
they do not reflect its real behaviour. Several authors (Carpinteri and
Mainardi, 1997; Hilfer, 2000; Gorenflo and Mainardi, 1998; Metzler and
Klafter, 2000) apply fractional calculus to the modelling of this type of
diffusion. This means that time and spatial derivatives in the classical
diffusion equation are replaced by fractional ones. In comparison to
derivatives of the integer order, which depend on the local behaviour of the
function, the derivatives of the fractional order accumulate the whole history of
this function.

In our previeus works (Ciesielski and Leszczynski, 2003, 2005), we presented
a solution to a~partial differential equation of the fractional order 
with the time fractional operator and the space ordinary operator,
respectively. Those solutions were based on the Finite Difference
Method (FDM) and are called the Fractional FDM (FFDM).

\section{Mathematical background}

\setcounter{equation}{0} 
In this paper we consider an ordinary differential equation of fractional 
order in the following form 
\begin{equation}
\frac{d ^{\alpha }}{d \left\vert x\right\vert _{\theta
}^{\alpha }}T(x)=0 \ \ \ \ \ x\in \mathbb{R}  \label{main_eqn}
\end{equation}%
where $T(x)$ is a field variable (i.e. field temperature), $\left(
d^{\alpha }/d \left\vert x\right\vert _{\theta }^{\alpha
}\right) T(x)$ is the Riesz-Feller fractional operator~(Metzer and Klafter,
2000; Samko et al, 1993), $\alpha $ is the real order of this operator and 
$\theta $ is the skewness parameter. According to (Gorenflo and Mainardi,
1998), the Riesz-Feller fractional operator is defined as
\begin{equation}
\frac{d^{\alpha }}{d \left\vert x\right\vert _{\theta
}^{\alpha }}T(x)={}_{x}D_{\theta }^{\alpha }T\left( x\right) ={}-\left[
c_{L}\left( \alpha ,\theta \right) \,_{-\infty }D_{x}^{\alpha }T\left(
x\right) +c_{R}\left( \alpha ,\theta \right) \,_{x}D_{+\infty }^{\alpha
}T\left( x\right) \right]  \label{def_RF}
\end{equation}%
for $0<\alpha \leq 2$, $\alpha \neq 1$ where 
\begin{eqnarray}
{}_{-\infty }D_{x}^{\alpha }T\left( x\right) &=&
\left( \dfrac{d}{dx} \right)^m \left[ _{-\infty }I_{x}^{m-\alpha }T\left( x\right) \right]
\\
{}_{x}D_{+\infty }^{\alpha }T\left( x\right) &=&
\left( -1\right)^m \left( \dfrac{d}{dx} \right)^m 
 \left[ _{x}I_{+\infty }^{m-\alpha }T\left( x\right) \right]   \notag
\end{eqnarray}
for $m\in \mathbb{N}$, $m-1<\alpha \leq m$, and the coefficients 
$c_{L}\left( \alpha ,\theta \right) $, $c_{R}\left( \alpha
,\theta \right) $ (for $0<\alpha \leq 2$, $\alpha \neq 1$, $%
\left\vert \theta \right\vert \leq \min \left( \alpha ,2-\alpha \right) $)
are defined as 
\begin{equation}
c_{L}\left( \alpha ,\theta \right) =\frac{\sin \dfrac{\left( \alpha -\theta
\right) \pi }{2}}{\sin \left( \alpha \pi \right) }\quad \quad c_{R}\left(
\alpha ,\theta \right) =\frac{\sin \dfrac{\left( \alpha +\theta \right) \pi 
}{2}}{\sin \left( \alpha \pi \right) }  \label{coeff_c}
\end{equation}

The fractional integral operators of the order $\alpha $: $_{-\infty
}I_{x}^{\alpha }T\left( x\right) $ and $_{x}I_{\infty }^{\alpha }T\left(
x\right) $ are defined as the left- and right-hand of Weyl's fractional
integrals (Carpinteri and Mainardi, 1997; Oldham and Spaner, 1974; Podlubny,
1999; Samko et al, 1993) whose definitions are 
\begin{equation}
_{-\infty }I_{x}^{\alpha }T\left( x\right) =\frac{1}{\Gamma (\alpha )}%
\int_{-\infty }^{x}\frac{T\left( \xi \right) }{(x-\xi )^{1-\alpha }}d\xi 
\end{equation}%
\begin{equation}
_{x}I_{\infty }^{\alpha }T\left( x\right) =\frac{1}{\Gamma (\alpha )}%
\int_{x}^{\infty }\frac{T\left( \xi \right) }{(\xi -x)^{1-\alpha }}d\xi \notag
\end{equation}

Considering Eq.~(2.1), we obtain the steady state of the classical diffusion 
equation for $\alpha =2$, i.e. the heat transfer equation. If $\alpha =1$, and the
parameter of skewness $\theta $ admits extreme values in~(\ref{coeff_c}),
the steady state of a transport equation is noted. Therefore we assume variations 
of the parameter $\alpha $ within the range $0<\alpha \leq 2$. Analysing behaviour
of the parameter $\alpha <2$ in Eq.~(2.1), we found some combination between
transport and propagation processes in steady states.

In this work we will consider equation (\ref{main_eqn}) in one dimensional
domain $\Omega : L\leq x\leq R$ with the boundary-value conditions of the
first kind (Dirichlet conditions) as 
\begin{equation}
\left\{ 
\begin{array}{ll}
x=L: & T\left( L\right) =g_{L} \\ 
x=R: & T\left( R\right) =g_{R}
\end{array}
\right.  
\label{bound_cond}
\end{equation}

\section{Numerical method}

\setcounter{equation}{0} 
According to the finite difference method (Hoffman, 1992; Majchrzak and
Mochnacki, 1996), we consider a~discrete from of equation
(\ref{main_eqn}) in space. The problem of solving equation (\ref{main_eqn})
lies in a properly approximation of Riesz-Feller derivative~(\ref{def_RF}) 
by a~numerical scheme.

\subsection{Approximation of the Riesz-Feller derivative}

In order to develop a discrete form of operator~(\ref{def_RF}), we
introduce a homogenous spatial grid $-\infty <\ldots
<x_{i-2}<x_{i-1}<x_{i}<x_{i+1}<x_{i+2}<\ldots <\infty $ with the step 
$h=x_{k}-x_{k-1}$. We denote a~value of the function $T$ in the point $x_{k}$ as
$T_{k}=T\left( x_{k}\right) $, for $k\in \mathbb{Z}$. For numerical
integration scheme, we assumed the trapezoidal rule and we used various
weighted numerical differential schemes for the first and second derivatives,
respectively.
The method of determination of the discrete form of this operator was described in
detail in work by Ciesielski (2005) and its final form is the following
\begin{equation}
_{x_{i}}D_{\theta }^{\alpha }T_{i}\approx \dfrac{1}{h^{\alpha }}%
\sum\limits_{k=-\infty }^{\infty }T_{i+k}w_{k}^{\left( \alpha ,\theta
\right) }  \label{discr_form}
\end{equation}%
where coefficients $w_{k}^{\left( \alpha ,\theta \right) }$, for $0<\alpha <1$,
have the form
\begin{eqnarray}
&&w_{k}^{\left( \alpha ,\theta \right) }=\dfrac{-1}{2\Gamma \left( 2-\alpha
\right) }  \cdot  \label{il_roz_Feller_1_w} \\
&& \cdot \left\{ 
\begin{array}{ll}
\left( \left( \left\vert k\right\vert +2\right) ^{1-\alpha }\lambda
_{1}+\left( \left\vert k\right\vert +1\right) ^{1-\alpha }\left( 2-3\lambda
_{1}\right) \right.  &  \\ 
\left. \quad +\left\vert k\right\vert ^{1-\alpha }\left( 3\lambda
_{1}-4\right) +\left( \left\vert k\right\vert -1\right) ^{1-\alpha }\left(
2-\lambda _{1}\right) \right) c_{L} & \text{ for }k\leq -2 \\ 
\left( 3^{1-\alpha }\lambda _{1}+2^{1-\alpha }\left( 2-3\lambda _{1}\right)
+3\lambda _{1}-4\right) c_{L}+\lambda _{1}c_{R} & \text{ for }k=-1 \\ 
\left( 2^{1-\alpha }\lambda _{1}-3\lambda _{1}+2\right) \left(
c_{L}+c_{R}\right) & \text{ for }k=0 \\ 
\left( 3^{1-\alpha }\lambda _{1}+2^{1-\alpha }\left( 2-3\lambda _{1}\right)
+3\lambda _{1}-4\right) c_{R}+\lambda _{1}c_{L} & \text{ for }k=1 \\ 
\left( \left( k+2\right) ^{1-\alpha }\lambda _{1}+\left( k+1\right)
^{1-\alpha }\left( 2-3\lambda _{1}\right) \right.  &  \\ 
\left. \quad  +k^{1-\alpha }\left( 3\lambda _{1}-4\right) +\left(
k-1\right) ^{1-\alpha }\left( 2-\lambda _{1}\right) \right) c_{R} & \text{
for }k\geq 2%
\end{array}%
\right.   \notag
\end{eqnarray}%
and for $1<\alpha \leq 2$, we obtain
\begin{align}
& w_{k}^{\left( \alpha ,\theta \right) }=\dfrac{-1}{2\Gamma \left( 3-\alpha
\right) } \cdot  \label{il_roz_Feller_2_w} \\
& \cdot \left\{ 
\begin{array}{ll}
\left( \left( \left\vert k\right\vert +2\right) ^{2-\alpha }\left( 2-\lambda
_{2}\right) +\left( \left\vert k\right\vert +1\right) ^{2-\alpha }\left(
4\lambda _{2}-6\right) + \right. & \\
\quad + \left\vert k\right\vert ^{2-\alpha }\left(6-6\lambda _{2}\right)+  
 \left( \left\vert k\right\vert -1\right) ^{2-\alpha
}\left( 4\lambda _{2}-2\right) + & \\ \quad + \left. \left( \left\vert k\right\vert -2\right)
^{2-\alpha }\left( -\lambda _{2}\right) \right) c_{L} & 
\text{for }k\leq -2 \\ 
\left( 3^{2-\alpha }\left( 2-\lambda _{2}\right) +2^{2-\alpha }\left(
4\lambda _{2}-6\right) -6\lambda _{2}+6\right) c_{L}+& \\ \quad +\left( 2-\lambda
_{2}\right) c_{R} & \text{for }k=-1 \\ 
\left( 2^{2-\alpha }\left( 2-\lambda _{2}\right) +4\lambda _{2}-6\right)
\left( c_{L}+c_{R}\right)  & \text{for }k=0 \\ 
\left( 3^{2-\alpha }\left( 2-\lambda _{2}\right) +2^{2-\alpha }\left(
4\lambda _{2}-6\right) -6\lambda _{2}+6\right) c_{R}+& \\ \quad + \left( 2-\lambda
_{2}\right) c_{L} & \text{for }k=1 \\ 
\left( \left( k+2\right) ^{2-\alpha }\left( 2-\lambda _{2}\right) +\left(
k+1\right) ^{2-\alpha }\left( 4\lambda _{2}-6\right) + \right. & \\
\quad + k^{2-\alpha }\left(6-6\lambda _{2}\right) 
+\left( k-1\right) ^{2-\alpha }\left( 4\lambda
_{2}-2\right) + & \\ \quad + \left. \left( k-2\right) ^{2-\alpha }\left( -\lambda _{2}\right)
\right) c_{R} & \text{for }k\geq 2
\end{array}%
\right.  \notag
\end{align}

Assuming $\alpha =2$ and $\theta =0$ and $c_{L}\left( 2,0\right) =
c_{R}\left( 2,0\right) = -1/2$, we obtain 
\begin{equation}
w_{k}^{\left( 2,0\right) }=\left\{ 
\begin{array}{ll}
0 & \text{for }k\leq -2 \\ 
1 & \text{for }k=-1 \\ 
-2 & \text{for }k=0 \\ 
1 & \text{for }k=1 \\ 
0 & \text{for }k\geq 2%
\end{array}%
\right. 
\end{equation}%
These coeeficients are identical as in the well known central difference
scheme for the second derivative.

The authors did not find exact values of the approximating coefficients
in literature. When $\alpha =1$, the Riesz-Feller operator is singular,
hence the problem occurs.

\subsection{Fractional FDM}
Having the discretization of the Riesz-Feller derivative in space done,
in this subsection we describe the finite difference method for 
equation (\ref{main_eqn}). Here we restrict the solution to
one dimensional space in comparison with the standard diffusion
equation where the discretization of the second derivative in space can
approximate the central difference of yhe second order. The differences 
appear in the setting of boundary conditions.

In the scheme of FDM, we replaced equation~(\ref{main_eqn}) by the
following formula 
\begin{equation}
\dfrac{1}{h^{\alpha }}\sum\limits_{k=-\infty }^{\infty }T_{i+k}w_{k}^{\left(
\alpha ,\theta \right) }=0  \label{diff_schem}
\end{equation}%
Here, for the unbounded domain we are obligated to solve a~system of algebraic
equations with an infinite dimension.

\subsubsection{Boundary value problem}

Present numerical scheme~(\ref{diff_schem}) with the included unbounded
domain $-\infty <x<\infty $ has no practical implementations in computer
simulations.

\begin{figure}[h]
\begin{center}
\includegraphics[width=0.6\textwidth]{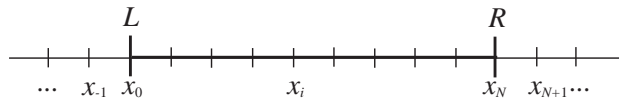}
\end{center}
\vspace{-0.5cm}
\caption{Spatial distribution of grid nodes}
\label{fig1}
\end{figure}

Now we show solution to this problem on the finite domain
$\Omega : L\leq x\leq R$ with boundary conditions~(\ref{bound_cond}). We divide
domain $\Omega $ into $N$ sub-domains with $h=(R-L)/N$. Figure~1 shows
a modified spatial grid. Here we can observe additional
'virtual' points in the grid placed outside of the domain $\Omega $. In
order to introduce the Dirichlet boundary conditions, we proposed a treatment
which is based on the assumption that values of the function $T$ in outside points
are identical as values in the boundary nodes $x_{0}$ or $x_{N}$ 
\begin{equation}
T\left( x_{k}\right) =\left\{ 
\begin{array}{ll}
T\left( x_{0}\right) =g_{L} & \text{for }k<0 \\ 
T\left( x_{N}\right) =g_{R} & \text{for }k>N
\end{array}%
\right. 
\end{equation}%
On the basis of above considerations, we modify expression~(\ref{discr_form})
for the discretization of the Riesz-Feller derivative. Thus we have 
\begin{equation}
_{x_{i}}D_{\theta }^{\alpha }T\left( x_{i}\right) \approx \dfrac{1}{%
h^{\alpha }}\left( \sum\limits_{k=-i}^{N-i}T_{i+k}w_{k}^{\left( \alpha
,\theta \right) }+g_{L}{s_{L}}_{i}^{\left( \alpha ,\theta \right) }+g_{R}{%
s_{R}}_{N-i}^{\left( \alpha ,\theta \right) }\right)   \label{diff_schem_mod}
\end{equation}%
for $i=1,\ldots ,N-1$, where 
\begin{eqnarray}
{s_{L}}_{j}^{\left( \alpha ,\theta \right) } &=&\sum\limits_{k=-\infty
}^{-j-1}w_{k}^{\left( \alpha ,\theta \right) }=\left\{ 
\begin{array}{ll}
c_{L}\left( \alpha ,\theta \right) \cdot \overset{\sim }{r}_{j} & \text{for 
}0<\alpha <1 \\ 
c_{L}\left( \alpha ,\theta \right) \cdot \overset{\approx }{r}_{j} & \text{for
}1<\alpha \leq 2%
\end{array}%
\right.   \notag \\ 
{} \\
{s_{R}}_{j}^{\left( \alpha ,\theta \right) } &=&\sum\limits_{k=j+1}^{\infty
}w_{k}^{\left( \alpha ,\theta \right) }=\left\{ 
\begin{array}{ll}
c_{R}\left( \alpha ,\theta \right) \cdot \overset{\sim }{r}_{j} & \text{for 
}0<\alpha <1 \\ 
c_{R}\left( \alpha ,\theta \right) \cdot \overset{\approx }{r}_{j} & \text{%
for }1<\alpha \leq 2
\end{array}%
\right.   \notag
\end{eqnarray}%
and 
\begin{eqnarray}
\overset{\sim }{r}_{j} &=&\dfrac{\left( j+2\right) ^{1-\alpha }\lambda
_{1}+\left( j+1\right) ^{1-\alpha }\left( 2-2\lambda _{1}\right)
+j^{1-\alpha }\left( \lambda _{1}-2\right) }{2\Gamma \left( 2-\alpha \right) 
} \notag \\ {} \\
\overset{\approx }{r}_{j} &=&\dfrac{\left( j+2\right) ^{2-\alpha }\left(
2-\lambda _{2}\right) +\left( j+1\right) ^{2-\alpha }\left( 3\lambda
_{2}-4\right) +j^{2-\alpha }\left( 2-3\lambda _{2}\right)}
{2\Gamma \left( 3-\alpha \right) }+\notag \\
&&+ \dfrac{\left( j-1\right)^{2-\alpha }\lambda _{2}}
{2\Gamma \left( 3-\alpha \right) } \notag 
\end{eqnarray}

Putting expression (\ref{diff_schem_mod}) to equation (\ref{main_eqn}),
we obtain the following finite difference scheme
\begin{eqnarray}
&&\sum\limits_{k=-i}^{N-i}T_{i+k}w_{k}^{\left( \alpha ,\theta \right) }+g_{L}%
{s_{L}}_{i}^{\left( \alpha ,\theta \right) }+g_{R}{s_{R}}_{N-i}^{\left(
\alpha ,\theta \right) }=0\text{, for }i=0,\ldots ,N  \notag \\
&&T_{0}=g_{L}  \label{schem_mod} \\
&&T_{N}=g_{R}  \notag
\end{eqnarray}%
The above scheme described by expression~(\ref{schem_mod}) can be written in
a~matrix form as 
\begin{equation}
\mathbf{A}\cdot \mathbf{T}=\mathbf{B}
\end{equation}%
where 
\begin{eqnarray}
\mathbf{A} &=&\left[ 
\begin{array}{cccccccc}
1 & 0 & 0 & 0 & \ldots & 0 & 0 & 0 \\ 
a_{-1} & a_{0} & a_{1} & a_{2} & \ldots & a_{N-3} & a_{N-2} & a_{N-1} \\ 
a_{-2} & a_{-1} & a_{0} & a_{1} & \ldots & a_{N-4} & a_{N-3} & a_{N-2} \\ 
a_{-3} & a_{-2} & a_{-1} & a_{0} & \ldots & a_{N-5} & a_{N-4} & a_{N-2} \\ 
a_{-4} & a_{-3} & a_{-2} & a_{-1} & \ldots & a_{N-6} & a_{N-3} & a_{N-4} \\ 
\vdots & \vdots & \vdots & \vdots & \ddots & \vdots & \vdots & \vdots \\ 
a_{-N+2} & a_{-N+3} & a_{-N+4} & a_{-N+5} & \ldots & a_{0} & a_{1} & a_{2}
\\ 
a_{-N+1} & a_{-N+2} & a_{-N+3} & a_{-N+4} & \ldots & a_{-1} & a_{0} & a_{1}
\\ 
0 & 0 & 0 & 0 & \ldots & 0 & 0 & 1%
\end{array}%
\right] \qquad \notag \\ {} \\
\mathbf{B} &=&\left[ g_{L}, b_{1}, b_{2}, b_{3}, b_{4}, \hdots, b_{N-2}, b_{N-1}
g_{R} \right]^{\mathrm{T}} \notag
\end{eqnarray}%
with 
\begin{eqnarray}
{ } & a_{j} = w_{j}^{\left( \alpha ,\theta \right) } \quad \quad \quad \quad \quad \quad & \text{ for }%
j=-N+1,\ldots ,N-1 \notag \\ {} \\
{ } & b_{j} = g_{L}{s_{L}}_{j}^{\left( \alpha ,\theta \right) }+g_{R}{s_{R}}%
_{j}^{\left( \alpha ,\theta \right) } & \text{ for }j=1,\ldots ,N-1 \notag
\end{eqnarray}%
and $\mathbf{T}=\left[ T_{0},T_{1},T_{2},\ldots ,T_{N}\right] ^{\mathrm{T}}$ is the
vector of unknown values of the function $T$.

We can observe that boundary conditions influence to all values of the
function in every node. In opposite to the second derivative over space,
which is approximated locally, the characteristic feature of the Riesz-Feller
and other fractional derivatives is the dependence on values of all domain
points. For $\alpha =2$ and $\theta =0$, our scheme is identical as with the well
known and used central difference scheme in space (Hoffman, 1992;
Majchrzak and Mochnacki, 1996).

The skewness parameter $\theta $ has significant influence on the
solution. For $\alpha \rightarrow 1^{+}$ and $\theta \rightarrow \pm 1^{+}$
one can obtain the classical ordinary differential equation.

\section{ Simulation results}

\setcounter{equation}{0} 
In this section, we present results of calculation.
In all presented simulations we assumed $0 \le x \le 1$. Figure~2 shows temperature 
profiles over space with boundary conditions $g_{L}=2$, $g_{R}=1$ for different
values $\alpha =\left\{ 0.1,0.5.0.75,1.01,1.25,1.5,1.75,2\right\} $ and
$\theta =0$. Figure~3 presents another example of the solution in 
which we assumed $\alpha =1.01$ and different values of the skewness parameter
$\theta =\left\{0,0.1,0.3,0.5,0.7,0.9,0.99\right\} $. 
\begin{figure}[!h]
\begin{center}
\includegraphics[width=0.7\textwidth]{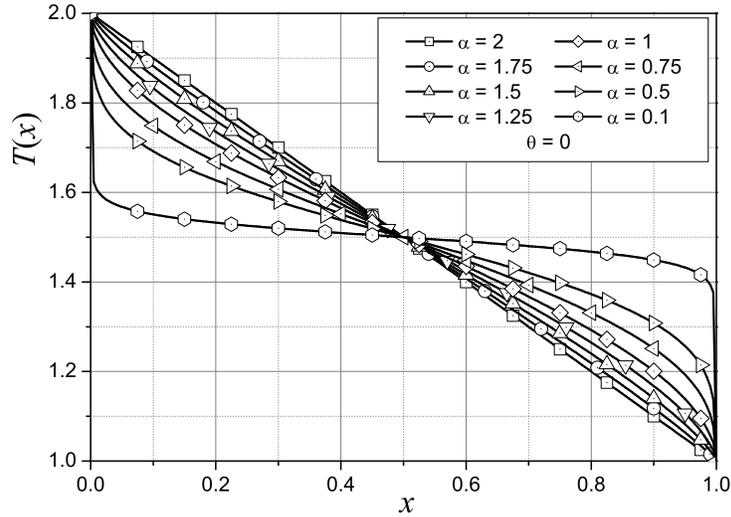}
\end{center}
\vspace{-0.5cm}
\caption{Spatial solution to equation (2.1) for different
       values of $\protect\alpha $ and $\protect\theta = 0$}
\label{fig2}
\end{figure}
\begin{figure}[!h]
\begin{center}
\includegraphics[width=0.7\textwidth]{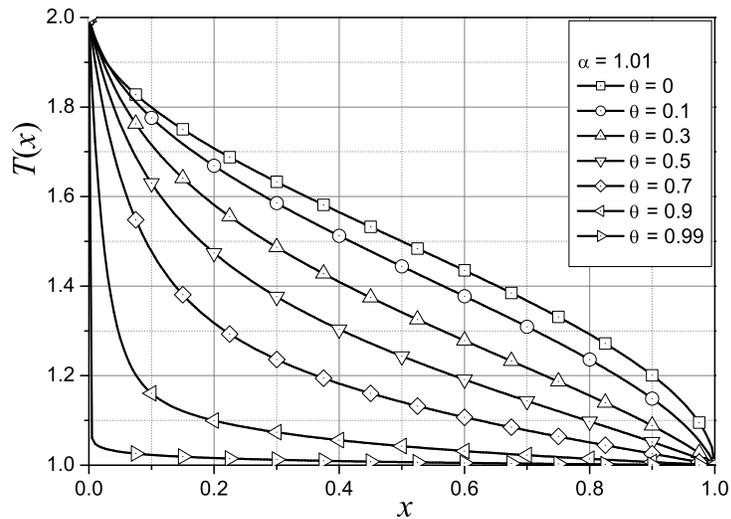}
\end{center}
\vspace{-0.5cm}
\caption{Spatial solution to equation (2.1) for a steady value of
       the parameter $\protect\alpha =1.01$ and different values of $\protect\theta$}
\label{fig3}
\end{figure}

The last example illustrates the heat transport in nanotubes. Figure~4 presents 
experimental data prerformed by Zhang and Li (2005) and results of
our numerical solution of equation (\ref{main_eqn}). For such a comparison we
assumed $\alpha=0.35$ and $\theta =-0.055$ to best fit the experimental data. It
should be noted that the temperature profile inside the nanotube deviates from
the profile obtained by solving the standard heat transfer equation.
\begin{figure}[h]
\begin{center}
\includegraphics[width=0.7\textwidth]{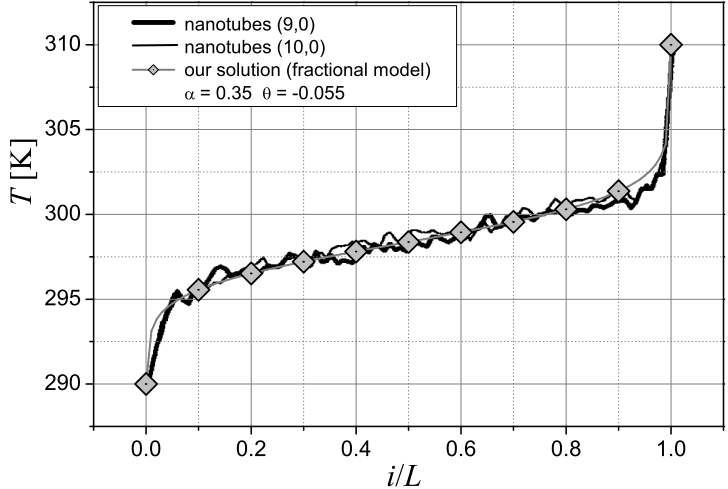}
\end{center}
\vspace{-0.5cm}
\caption{Comparison of temperature profiles of nanotubes measured by Zhang
         and Li (2005) and obtained by numerical solution to equation (2.1)}
\label{fig4}
\end{figure}

\section{ Conclusions}

Summing up, we proposed a fractional finite difference method 
for the fractional diffusion equation with the Riesz-Feller 
fractional derivative which is an extension of the standard 
diffusion. We analysed a linear case of the steady state of
the anomalous diffusion equation, and in the future we will 
work on non-linear cases. We obtained the FDM scheme which 
generalises classical scheme of FDM for the diffusion equation. 
Moreover, for $\alpha = 2$, our solution is equivalent with 
that obtained by the classical finite difference method.

Analysing plots included in this work, we can see that in the
case $\alpha = 2$, solution of Eq. (2.1) is a linear function.
In the other cases, when $\alpha < 2$, solutions are non-linear 
functions. When we analyse the probability density function 
generated by fractional diffusion equation for $\alpha < 2$,
we observe a long tail of distribution. In this way we can 
simulate same rare and extreme events which are characterised 
by arbitrary very large values of particle jumps.

Analysing the changes in the skewness parameter $\theta$, 
we observed interesting behaviour in the solution. For 
$\alpha \rightarrow 1^{+}$ and for $\theta \rightarrow \pm 1^{+}$,
we obtained the steady state of the first order wave equation.
For $\theta \in (0,1)$ (with restrictions to the order $\alpha$),
we generated non-symmetric solutions. It should be noted that 
we can good approximate the temperature profile inside
the nanotubes using solution of Eq. (2.1).

\end{document}